\documentclass[12pt,reqno]{amsart}
\usepackage{amsmath,amssymb,graphicx,amscd}
\usepackage{mathrsfs}
\usepackage{enumerate}
\usepackage{mathrsfs}
\usepackage{dsfont}
\usepackage{amsfonts}
\usepackage{fullpage}
\usepackage{color}
\renewcommand{\P}{\mathbf{P}}
\newcommand{\Q}{\mathbf{Q}}
\newcommand{\E}{\mathbf{E}}
\newcommand{\I}{\mathbf{1}}
\newcommand{\FF}{\mathbb{F}}
\newcommand{\GG}{\mathbb{G}}
\newcommand{\F}{\mathcal{F}}
\newcommand{\G}{\mathcal{G}}

\newcommand{\finproof}{$\square$}

\newcommand{\dd}{{\mathrm{d}}}            %differential

\newtheorem{thm}{Theorem}[section]
\newtheorem{cor}[thm]{Corollary}
\newtheorem{lem}[thm]{Lemma}
\newtheorem{prop}[thm]{Proposition}
\theoremstyle{definition}
\newtheorem{defn}[thm]{Definition}
\theoremstyle{remark}
\newtheorem*{rem}{Remark}

\theoremstyle{remark}
\newtheorem{ex}[thm]{Example}
\numberwithin{equation}{section}

\numberwithin{equation}{section}

\begin{document}

\title{Default times, non arbitrage conditions and change of probability measures}

\author{Delia Coculescu}
\address{ETHZ \\ Departement Mathematik, R\"{a}mistrasse 101\\ Z\"{u}rich 8092, Switzerland.}
\email{delia.coculescu@math.ethz.ch}

\author{Monique Jeanblanc}
\address{Universit\'e  d'Evry Val d'Essonne \\ D\'epartement de Math\'ematiques  \\
rue Jarlan F-91025 Evry Cedex France.\thanks{This research
benefited from the support of the ``Chaire Risque de Cr\'edit'',
F\'ed\'eration Bancaire Fran\c caise.}}
\email{monique.jeanblanc@univ-evry.fr}

\author{Ashkan Nikeghbali}
\address{Institut f\"ur Mathematik,
Universit\"at Z\"urich, Winterthurerstrasse 190,
CH-8057 Z\"urich,
Switzerland}
\email{ashkan.nikeghbali@math.uzh.ch}

\subjclass[2000]{15A52} \keywords{Default modeling, credit risk models, random times, enlargements of filtrations, immersed filtrations, no arbitrage conditions, equivalent change of measure.}

\begin{abstract}
In this paper we give  a financial justification, based on non
arbitrage conditions, of the $(H)$ hypothesis in default
time modelling. We also show how the $(H)$ hypothesis is affected
by an equivalent change of probability measure.  The main
technique used here is the theory of progressive enlargements of
filtrations.
\end{abstract}

\maketitle

\section{Introduction}

In this paper we study the stability of the $(H)$ hypothesis (or immersion property) under equivalent changes of probability measures. Given two filtrations $\FF \subset \GG$, we say that $\FF$ is immersed in $\GG$ if all $\FF$-local martingales are $\GG$-local martingales. In the default risk literature, the filtration $\GG$ is obtained by the progressive enlargement of $\FF$ with a random time (the default time) and the immersion property under a risk-neutral measure appears to be a suitable non arbitrage condition (see \cite{bsjeanblanc} and \cite{jealec08}). Because in general immersion is not preserved under equivalent changes of probability measures (see \cite{yorjeulinens} and \cite{emery}), the reduced form default models are usually specified directly under a given risk-neutral measure. 

However,  it seems crucial to understand how the immersion property is modified under an equivalent change of probability measure. This is important not only because the credit markets are highly incomplete, but also because the physical default probability  as well appears to play an important role in presence of incomplete information, as emphasized by a more recent body of literature, initiated by \cite{DuffLand} (see also \cite{GieseGold1}, \cite{JeanblValch04}, \cite{cocgemjea08}, \cite{FreySchmidt} among others). This imperfect information modeling approach proposes to rely on accounting information, and to incorporate the imperfect information about the accounting indicators, in computing the credit spreads. The default intensities are computed endogenously, using the available observations about the firm.  Some of the constructions do not satisfy the immersion property (\cite{Kusuok}, \cite{GuoJarrZeng}). It is therefore important to understand the role of the immersion property for pricing.

More generally, our goal in this paper is to provide efficient and precise tools from martingale theory and the general theory of stochastic processes to model default times: we wish to justify on economic grounds the default models which use the technique of  progressive enlargements of filtrations, and to explain the reasons why such an approach is useful. We  provide and study (necessary and) sufficient conditions for a market model to be arbitrage free in presence of default risk. More precisely, the paper is organized as follows:

In section 2, we describe the financial framework which  uses the enlargements of filtrations techniques and introduce the corresponding non arbitrage conditions. In section 3, we present the useful tools form the theory of the progressive enlargements of filtrations.

Eventually, we study how the immersion property is affected under equivalent changes of probability measures. In section 4, we give a simple proof of the not well known fact (due to Jeulin and Yor \cite{yorjeulinens}) that immersion is preserved under a change of probability measure whose Radon-Nikod\'ym density is $\F_\infty$-measurable. Using this result, we show that a sufficient non arbitrage condition is that the immersion property should hold under an equivalent change of measure (not necessarily risk-neutral). Then, using a general representation property for $\GG$ martingales (section 5), we characterize the class of equivalent changes of probability measures which preserve the immersion property when the random time $\tau$ avoids the $\FF$ stopping times (section 6), thus extending the results of Jeulin and Yor \cite{yorjeulinens} in our setting.  Eventually,  we show how the Az\'ema supermartingale is computed for a large class of equivalent changes of measures.

\section{Non arbitrage conditions}\label{sec::finmodel}

In this section we briefly comment some non arbitrage conditions appearing in the default models that use the progressive enlargement of a reference filtration (for further discussion in complete default-free markets see \cite{bsjeanblanc} and \cite{jealec08}). All the notions from the theory of enlargements of filtrations we use in this section are gathered in the next section.
In default modeling, the technique of progressive enlargements of filtrations has been introduced by Kusuoka \cite{Kusuok} and further developed in Elliott, Jeanblanc and Yor \cite{elliotjeanbyor}. It consists in a two step construction of the market model, as follows.

Let $(\Omega, \G, \FF=(\F_t)_{t\geq 0}, \P) $ be a filtered probability space satisfying the usual hypothesis. For us, the probability $\P$ stands for the physical measure under which financial events and prices are observed. Let   $\tau$ be a random time:  it is a $\G$-measurable random variable which usually represents the default time of the company. It is not an $\FF$ stopping-time. Let $\GG=(\G_t)_{t\geq 0}$ be the filtration obtained by progressively enlarging the filtration $\FF$ with the random time $\tau$. Obviously, $\forall t\geq 0$, $\F_t\subset\G_t\subset\G$.

Usually,  the filtration $\GG$ plays the role of the market filtration (and is sometimes called the full market
filtration), meaning that the price processes are $\GG$-adapted, and the pricing is performed with respect to this filtration. On the other hand, the definition of the filtration $\FF$ (called the reference filtration) is not always clear in the literature so far, and several interpretations can be given.

Let us now suppose that the reference filtration $\FF$ contains the \textit{market price} information which an investor is using for evaluating some defaultable claims. Typically, this is the natural filtration of a vector of semi-martingales  $S=(S_t)_{t \geq 0}$, with $S:=(S^1,..S^n)$. The vector $S$ is recording the prices history of observable default-free (with respect to $\tau$) assets which are sufficiently liquid to be used for calibrating the model. Here, we may include assets without default risk, as well as assets with a different default time than $\tau$, typically  assets issued by other companies than the one we are analyzing. We shall call  $\tau$-default-free assets the components of $S$ , since these are not necessarily assets without default risk.

As usual, we let $S^0$ stand for the locally risk-free asset (i.e., the money market account); the remaining assets are risky. We denote by $\Theta(\FF,\P)$ the set of all equivalent local martingale measures for the num\'eraire $S^0$, i.e.:
\[
\Theta(\FF,\P)=\left\{\Q \sim \P \text{ on }\G| \frac{S}{S^0}=\left(
\frac{S^1}{S^0},..., \frac{S^{n}}{S^0}\right) \text{ is a
}(\FF,\Q)-\text{local martingale}\right\},
\]
and we will suppose that $\Theta(\FF,\P)$ is not empty in order to ensure absence of arbitrage opportunities (see \cite{delbscha94}). Notice that, because we shall work with different filtrations, we prefer to always define the probability measures on the sigma algebra $\G$. In this way, we avoid dealing with extensions of a probability measure. When the $\FF$-market is complete, all the measures belonging to  $\Theta(\FF,\P)$ have the same $\FF$ restriction.

In practice, investors might use different
information sets than $\FF$, say $\GG$. In this case, they can
construct $\GG$-portfolios and $\GG$-strategies. Then, from the
viewpoint of the arbitrage theory, one needs to understand what
the relevant prices become in a different filtration.

In particular, some investors may use more than the information in $\FF$ for constructing the portfolios. For instance, they might take into account the macro-economic environment, or firm specific
accounting information which is not directly seen in the prices.
In this case $\FF\subset\GG$. Denote:
\[
\Theta(\GG,\P)=\left\{\Q \sim \P \text{ on }\G| \frac{S}{S^0}=\left( \frac{S^1}{S^0},..., \frac{S^{n}}{S^0}\right) \text{ is a }(\GG,\Q)-\text{local martingale}\right\}.
\]
Are there (local) martingale measures for the $\GG$-informed traders? One has to understand what the $\FF$-martingales become in a larger filtration. There is not a general answer to this question: in general martingales of a given filtration are not semi-martingales in a larger filtration (\cite{yorjeulin}). However, from a purely economic point of view, if one assumes that the information contained in $\GG$ is available for all the investors without cost (i.e., this is public information), then the non arbitrage condition becomes: $$\Theta(\GG,\P)\neq \emptyset.$$ This is coherent with the semi-strong form of the market efficiency, which says that a price process fully reflects all relevant information that is publicly available to investors. This means that publicly available information cannot be used in order to obtain arbitrage profits.

Let us now come to the particular case of the default models, where  $\FF$ stands for the information about the prices of $\tau$-default-free assets. In general, $\tau$ is not an $\FF$ stopping time and for the purpose of pricing defaultable claims, the progressively enlarged filtration $\GG$ has to be introduced. As an illustration, let us take the filtering model introduced by Kusuoka:

\begin{ex}\label{ExKusuoka}
\textbf{Kusuoka's filtering model (1999):} Let $(B^1_t,B^2_t)_{t\in[0,T]}$ a 2-dimensional Brownian motion. The default event is triggered by the following process (for instance the cash flow balance of the firm, or assets' value):
\[
dX_t:=\sigma^1(t,X_t)dB^1_t+b(t,X_t)dt, \quad X_0=x_0.
\]
Let $\tau :=\inf\{t\in [0,T]| X_t=0\}$ be the default time.
Suppose that the market investors do not observe $X$, but instead the following process:
\[
dY_t:=\sigma^2(t,Y_t)dB^2_t+\mu(t,X_{t\wedge\tau},Y_t)dt, \quad Y_0=y_0.
\]
The process $Y$ might be a $\tau$-default-free asset price that is correlated with the defaultable assets value. For instance, suppose $X$ is the assets value of an oil company. Then, the oil price is an important piece of information to take into account when estimating the default risk of the company. Then $Y$ can be the spot price of oil. The reference filtration is $\F_t:=\sigma(Y_s,s \leq t)$ and the market filtration is constructed as $\G_t:=\F_t \vee \sigma(\tau \wedge s, s\leq t).$\end{ex}

As Kusuoka pointed out, the above example does not fulfill the immersion property. It is natural to investigate if such a model is arbitrage-free.

Let us assume that $\Theta(\FF,\P)$ is not empty, i.e., the $\tau$-default-free market is arbitrage free, and let us introduce the following alternative non arbitrage conditions:

\begin{enumerate}
\item There exists $\Q\in \Theta(\FF,\P)$ such that every $(\FF,\Q)$-martingale is a $(\GG,\Q)$-martingale, i.e., the immersion property holds under $\Q$.

\item There exists a measure $ \Q \sim \P$ such that every $(\FF,\P)$-martingale is a $(\GG,\Q)$-martingale.
\end{enumerate}

The idea behind both conditions is that, since default events are public information, an investor who uses this information to decide his trading strategy should not be able to make arbitrage profits. Condition (1) says that there is (at least) one martingale measure in common for an investor who uses information from default (filtration $\GG$) in his trading and a less informed one, who is only concerned with $\tau$-default-free prices levels when trading (filtration $\FF$). Condition (2) looks at first sight less restrictive, by only saying that for each such type of investor there exists a martingale measure (but which could a priori  be different).
A closer inspection tells us that the two conditions are in fact equivalent. This equivalence  will be proved in section 4 where we also show that these conditions are equivalent to:
\begin{enumerate}
\item[(3)] There exists $\Q\sim\P$ such that the immersion property holds under $\Q$.
\end{enumerate}
In other words, as soon as the immersion property holds under an equivalent probability measure, immersion holds as well under (at least) one $\FF$ risk neutral measure. Furthermore, $\Theta(\GG,\P)$ is not empty, i.e.,  non arbitrage holds for the defaultable market. Hence, the immersion property is an important non arbitrage condition to study.

Note also that the conditions listed above are sufficient for $\Theta(\GG,\P)$ to be not empty but not necessary. One only needs that the martingale invariance property holds for the price processes $S$, not for all $\FF$ local martingales. Thus, when the $\FF$ market is incomplete, weaker conditions can be stated. We now recall some important facts from the theory of progressive enlargements of filtrations which are relevant to our study.

\section{Basic facts about random times and  progressive enlargements of filtrations}\label{sec::tpsaleatoires}

In this section, we recall some important facts from the general theory of stochastic processes which we shall need in the sequel.
We assume we are given a filtered probability space $(\Omega,\G,\FF=( \mathcal{F}_{t})_{t \geq 0},\P)$  satisfying the usual assumptions
\begin{defn}
A random time $\tau$ is a nonnegative random variable $\tau:\left(\Omega,\mathcal{G}\right)\rightarrow[0,\infty]$.
\end{defn}

The theory of progressive enlargements of filtrations was
introduced to study properties of random times which are not
stopping times: it originated in a paper by Barlow \cite{barlow}
and was further developed by Yor and Jeulin, \cite{yorgrossis}
\cite{yorjeulin} \cite{jeulin79,jeulin}. For further  details, the
reader can also refer to \cite{jeulinyor} which is written in
French or to \cite{columbia} or \cite{protter} chapter VI for an
English text. This theory gives the decomposition of local
martingales of the initial filtration $\FF$ as semimartingales of
the progressively enlarged one $\GG$. More precisely, we enlarge
the initial filtration $\FF$ with the one generated by the
process $\left( \tau \wedge t\right) _{t\geq 0}$, so that the new
enlarged filtration $\GG=\left( \mathcal{G}_{t}\right) _{t\geq
0}$ is the smallest filtration (satisfying the usual
assumptions) containing $\FF$\ and making $\tau $\ a stopping
time, i.e.,
$$\G_{t}=\mathcal{K}_{t+}^{o} \quad \text{where} \quad \mathcal{K}_{t}^{o}=\F_{t}\vee\sigma(\tau\wedge t).$$

A few processes  play a crucial role in our discussion:
\begin{itemize}
\item the $\FF$ supermartingale
\begin{equation}
Z_t^\tau=\mathbf{P}\left[ \tau >t\mid \mathcal{F}_{t}\right]
\label{surmart}
\end{equation}
chosen to be c\`{a}dl\`{a}g, associated with $\tau $\ by Az\'{e}ma (\cite{azema});

\item the $\FF$ dual optional and predictable projections of the process $1_{\left\{ \tau \leq t\right\} }$,
denoted respectively by $A_{t}^{\tau }$ and $a_{t}^{\tau }$;

\item the c\`{a}dl\`{a}g martingale
\begin{equation*}
\mu _{t}^{\tau }=\mathbf{E}\left[ A_{\infty }^{\tau }\mid \mathcal{F}_{t}\right] =A_{t}^{\tau }+Z_{t}^{\tau }.
\end{equation*}

\item the Doob-Meyer decomposition of (\ref{surmart}):
\begin{equation*}
Z_t^\tau=m_{t}^{\tau }-a_{t}^{\tau },
\end{equation*} where $m^{\tau}$ is an $\FF$-martingale.
\end{itemize}

In the credit risk literature, the hazard process is very often used:
\begin{defn}
Let $\tau$ be  a random time such that $Z_t^\tau>0$, for all $t\geq0$ (in particular $\tau$ is not an $\FF$- stopping time). The nonnegative stochastic process $\left(\Gamma_t\right)_{t\geq0}$ defined by:
$$\Gamma_t=-\ln Z_t^\tau,$$
is called the \emph{hazard process}.
\end{defn}

It is important to know how the $\FF$ local martingales are
affected under the progressive enlargement of filtrations: in
general, for an arbitrary random time, an $\FF$ local martingale
is not a $ \GG$ semimartingale (see \cite{jeulin}, (\cite{yorjeulin}). However, we
have the following general result:

\begin{thm}[Jeulin-Yor \cite{yorjeulin}]\label{decotempsqcq}
Every $\FF$ local martingale $\left( M_{t}\right) $, stopped at $\tau $, is a $\GG$ semimartingale, with canonical decomposition:
\begin{equation}
M_{t\wedge \tau }=\widetilde{M}_{t}+\int_{0}^{t\wedge \tau
}\frac{d\left\langle M,\mu ^{\tau }\right\rangle_{s}}{Z_{s-}^{\tau }}  \label{decocanonique}
\end{equation}
where $( \widetilde{M}_{t}) $\ is a $\GG$ local martingale.
\end{thm}

Moreover, the Az\'ema supermartingale is the main tool for computing  the $\GG$ predictable compensator of $\I_{\tau\leq t}$:

\begin{thm}[\cite{yorjeulin}]\label{calccomp}
Let $H$ be a bounded $\GG$ predictable process. Then
$$H_{\tau}\I_{\{\tau\leq t\}}-\int_{0}^{t\wedge\tau}\dfrac{H_{s}}{Z_{s-}^{\tau}}da_{s}^{\tau}$$
is a $\GG$ martingale. In particular, taking $H\equiv1$, we find that:
$$N_t:=\I_{\{\tau\leq t\}}-\int_{0}^{t\wedge\tau}\frac{1}{Z_{s-}^{\tau}}da_{s}^{\tau}$$
is a $\GG$ martingale.
\end{thm}

The following assumptions are often encountered in the literature on enlargements of filtrations or the modelling of default times:

\begin{itemize}
\item The $\mathbf{(H)}$-hypothesis: every $\FF$ martingale is a $\GG$ martingale.
We say that the filtration $\FF$ is immersed in $\GG$.

\item Assumption $\mathbf{(A)}$: the random time $\tau $ \underline{a}voids every $\FF$ stopping time $T$, i.e. $\mathbf{P}\left[ \tau =T\right] =0$.
\end{itemize}
When one assumes that the random time $\tau$ avoids $\FF$ stopping times, then one further has:
\begin{lem}[\cite{yorjeulin}, \cite{jeulin}]\label{lem:evitement}
If $\tau $\ avoids  $\FF$ stopping times (i.e. condition
$\mathbf{(A)}$ is satisfied),  then $A^{\tau }=a ^{\tau }$  and
$A$ is continuous. The  $\GG$ dual predictable projection of the
process $\I_{\{\tau\leq t\}}$ is continuous and $\tau$ is a
totally inaccessible stopping time.
\end{lem}

We now state several useful equivalent characterizations of the $\mathbf{(H)}$ hypothesis in the next theorem. Note that except the last equivalence, the results are true for any filtrations $\FF$ and $\GG$ such that $\F_t\subset\G_t$.
The theorem is a combination of results by Br\'emaud and Yor \cite{bremaudyor} and also by  Dellacherie and Meyer \cite{delmeyfil} in the special case when the larger filtration is obtained by progressively enlarging the smaller one with a random time.

\begin{thm}[Dellacherie-Meyer \cite{delmeyfil} and Br\'{e}maud-Yor \cite{bremaudyor}]\label{lecasdregeneral} The following are equivalent:

\begin{enumerate}
\item $(H)$: every $\FF$ martingale is a $\GG$ martingale;

\item For all bounded $\F_{\infty }$-measurable random variables $\mathbf{F}$\ and all bounded $\mathcal{G}_{t}$-measurable random variables $\mathbf {G}_{t}$, we have
\begin{equation*}
\mathbf{E}\left[ \mathbf{FG}_{t}\mid \mathcal{F}_{t}\right]
=\mathbf{E}\left[ \mathbf{F}\mid \mathcal{F}_{t}\right]
\mathbf{E}\left[ \mathbf{G}_{t}\mid \mathcal{F}_{t}\right] .
\end{equation*}

\item For all bounded $\mathcal{F}_{\infty }$ measurable random variables $\mathbf{F}$,
\begin{equation*}
\mathbf{E}\left[ \mathbf{F}\mid \mathcal{G}_{t}\right] =\mathbf{E}\left[ \mathbf{F}\mid \mathcal{F}_{t}\right] .
\end{equation*}

\item For all $s\leq t$,
\begin{equation*}
\mathbf{P}\left[ \tau \leq s\mid \mathcal{F}_{t}\right]
=\mathbf{P}\left[ \tau \leq s\mid \mathcal{F}_{\infty }\right] .
\end{equation*}
\end{enumerate}
\end{thm}

We now indicate some consequences of the condition $\mathbf{(A)}$.

\begin{cor}
Suppose that the immersion property holds. Then
$Z^{\tau}=1-A^{\tau}$ is a decreasing process. Furthermore, if
$\tau$ avoids stopping times, then $Z^{\tau} $ is continuous.
\end{cor}
\proof This is an immediate consequence of Theorem \ref{lecasdregeneral} and Lemma \ref{lem:evitement}.
\finproof

\begin{rem}

\begin{itemize}
\item[(i)] It is known that if $\tau$ avoids $\FF$ stopping times, then $Z^\tau$ is continuous and decreasing if and only if $\tau$ is a pseudo-stopping time (see \cite{AshkanYor} and \cite{coculescunikeghbali}). 

\item[(ii)] When the immersion property holds and $\tau$  avoids the $\FF$ stopping times, we have from the above corollary and Theorem \ref{calccomp} that the $\GG$ dual predictable projection of $\I_{\{\tau\leq t\}}$ is $\log\left(\frac{1}{Z_{t\wedge\tau}^\tau}\right)$.
\end{itemize}
\end{rem}

\section{Immersion property and equivalent changes of probability measures: first results}

In the remainder of the paper, the setting is the one of the previous section:  $(\Omega,\G,\FF=( \mathcal{F}_{t})_{t \geq 0},\P)$  is a filtered probability space satisfying the usual assumptions, $\tau$ is a random time and $\GG=( \mathcal{G}_{t})_{t \geq 0}$ is the progressively enlarged filtration which makes $\tau$ a stopping time.

\textbf{Notations.} We note $\FF \overset{\P}{\hookrightarrow} \GG$ for $\FF$ is immersed in $\GG$ under the probability measure $\P$. Let $\mathcal{I}(\P)$ be the set of all probability measures $\Q$ which are equivalent to $\P$ and such that $\FF \overset{\Q}{\hookrightarrow} \GG$.

We would now like to see how the immersion property is affected by equivalent changes of probability measures. Let $\Q$ be a probability measure which is equivalent to
$\P$ on $\G$, with $\rho= d\Q/d\P$.   Define:
\begin{equation}\label{densities}
\dfrac{d\mathbf{Q}}{d\mathbf{P}}\Big|_{\mathcal{F}_{t}}=e_{t};\quad\dfrac{d\mathbf{Q}}{d\mathbf{P}}\Big|_{\mathcal{G}_{t}}=E_{t}.
\end{equation}
We shall always consider c\`{a}dl\`{a}g versions of the martingales $e$ and $E$.

What can one say about the $(\FF,\Q)$ martingales when considered in the filtration $\GG$? A simple application of Girsanov's theorem yields:
\begin{prop}
Assume that $\FF \overset{\P}{\hookrightarrow} \GG$. Let $\Q$ be a
probability measure which is equivalent to $\P$ on $\G$. Then
every $(\FF,\Q)$ semimartingale is a $(\GG,\Q)$ semimartingale.
\end{prop}

The decomposition of the $(\FF,\Q)$-martingales in the larger filtration can be found by applying twice Girsanov's theorem, respectively in the filtration $\FF$ and  then in the filtration $\GG$:

\begin{thm}[Jeulin-Yor \cite{yorjeulinens}]\label{JY_decomp}

Assume that $\FF \overset{\P}{\hookrightarrow} \GG$. With the notation introduced in (\ref{densities}), if $\left(X_{t}\right)$ is an $\left(\FF,\mathbf{Q}\right)$-local martingale, then the stochastic process:

$$I^{X}_t :=X_{t}+\int_{0}^{t}\dfrac{E_{s-}}{E_{s}}\left(\dfrac{1}{e_{s-}}d[X,e]_{s}-\dfrac{1}{E_{s-}}d[X,E]_{s}\right)$$
is a $\left(\GG,\Q\right)$-local martingale. Note that
$$I^{X}_t=X_{t}+\int_{0}^{t}\frac{1}{\eta_{s-}}d[X,\eta]_{s}$$
where $\eta=e/E$ is a $(\GG,\Q)$-martingale.
\end{thm}

The decomposition above depends on the ratio $\eta=e/E$, hence on the initial probability $\P$. Can one instead find a decomposition involving the $\Q$-Az\'ema supermartingale? To answer this question, one has to understand on the one hand, how the Az\'ema supermartingale is affected by equivalent changes of measure and on the other hand, what measures preserve the immersion property.

We now give as a consequence of the Theorem \ref{lecasdregeneral} an invariance property for the Az\'ema supermartingale associated with $\tau$ for a particular class of equivalent changes of measure:

\begin{prop}\label{invariancecomp}
Let $\FF \overset{\P}{\hookrightarrow} \GG$ and let  $\Q$ be a probability measure which is equivalent to $\P$ on $\G$.  If $d\Q/d\P$ is $\F_\infty$-measurable, then:
$$\Q(\tau>t|\F_t)=\P(\tau>t|\F_t)=Z^\tau_t.$$
Consequently, the predictable compensator of $\I_{\{\tau\leq t\}}$ is unchanged under such equivalent changes of probability measures, i.e.,
$$N_t=\I_{\{\tau\leq t\}}-\int_{0}^{t\wedge\tau}\frac{da_s^\tau}{Z_{s-}^{\tau}}$$ is a $\GG$-martingale under $\P$ and $\Q$. Moreover, $\FF \overset{\Q}{\hookrightarrow} \GG$.
\end{prop}

\proof We have, for $s\leq t$:
$$\Q(\tau>s|\F_t)
=\dfrac{\E[\rho\I{\tau>s}|\F_t]}{\E[\rho|\F_t]};$$ and from
Theorem \ref{lecasdregeneral} (2), we have:
$\E[\rho\I_{\tau>s}|\F_t]=\E[\rho|\F_t]\E[\I_{\tau>s}|\F_t]=\E[\rho|\F_t]\P(\tau>s|\F_t)$,
and hence
$\Q(\tau>s|\F_t)=\P(\tau>s|\F_t)=\P(\tau>s|\F_\infty)=\Q(\tau>s|\F_\infty)$.

The result then follows from Theorem \ref{lecasdregeneral} (4).
\finproof

Now, we are able to deduce from Theorem \ref{JY_decomp} the following equivalence:

\begin{prop}
We do not assume immersion under $\P$. The following conditions are equivalent:
\begin{enumerate}
\item $\mathcal{I}(\P)\neq\emptyset$.
\item There exists $\Q\sim\P$ such that every $(\FF,\P)$  martingale is a $(\GG,\Q)$  martingale.
\end{enumerate}
\end{prop}

\proof

$(1) => (2)$. Suppose $\tilde \Q \in \mathcal{I}(\P)$. We apply Theorem \ref{JY_decomp} but
(unfortunately!) with the role of $\P$ taken here by $\tilde \Q$: If $X$ is a $(\FF,\P)$ local martingale then $X_t+\int_{0}^{t}\frac{1}{\eta_{s-}}d[X,\eta]_{s}$ is a $(\GG,\P)$ local martingale. Since $\eta$ is a $(\GG,\P)$-martingale, one can define   $d \Q = \eta_t \cdot d \P$ on $\G_t$. Applying Girsanov's theorem again, we obtain that $X_t+\int_{0}^{t}\frac{1}{\eta_{s-}}d[X,\eta]_{s}-\int_{0}^{t}\frac{1}{\eta_{s-}}d[X,\eta]_{s}=X_t$ is a  $(\GG,\Q)$ local martingale, hence (2) holds.

$(2) => (1)$. Let $m$ be any $(\FF,\P)$  martingale, hence by the
statement (2), $m$ is also a $(\GG,\Q)$ martingale, which is
$\FF$-adapted. It follows that $m$ is also an $(\FF, \Q)$
martingale.
In particular, the $(\FF,\P)$ martingale $e_t=\frac{d\Q}{d\P}|\F_t
$ is also an $(\FF,\Q)$-martingale. From Girsanov's theorem, this
is possible if and only if $[e,e]=0$, which implies that $e=1$,
hence $d \Q =   d \P$ on $\F_t$. Hence all $(\FF,\Q)$-martingales are  $(\FF,\P)$-martingales, hence  $(\GG,\Q)$-martingales, i.e., (H) holds under $\Q$. \finproof

Let us now go back to the financial framework of Section \ref{sec::finmodel}, where $\P$ stands for the physical measure, and let us analyze the non arbitrage conditions introduced there. We suppose that $\Theta(\FF,\P)$ is not empty, i.e., the $\FF$-market is arbitrage free. Now, we show that if there exists an equivalent probability measure such that immersion holds, then there exists as well a risk neutral one such that immersion holds, in other words:

\begin{prop}
If $ \Theta(\FF,\P)$ and  $\mathcal{I}(\P)$ are not empty, then $\Theta(\FF,\P)\bigcap\mathcal{I}(\P)\neq\emptyset$.
\end{prop}
\proof
Suppose $\Q\in\mathcal{I}(\P)$ and $\P^1\in \Theta(\FF,\P)$ such that $\P^1\notin\mathcal{I}(\P)$. Denote $d\P^1/d\Q|_{\F_\infty}=A$ and introduce $\P^2$ as $d\P^2/d\Q=A$. Since $A$ is $\F_{\infty}$-measurable, by Proposition \ref{invariancecomp}, $\P^2\in\mathcal{I}(\P)$. Moreover $\P^2\in \Theta(\FF,\P)$  since $d\P^2/d\P^1|_{\F_\infty}=1$. \finproof

The two above propositions tell us that a sufficient non arbitrage condition for the financial market introduced in Section \ref {sec::finmodel} is:   $\mathcal{I}(\P)\neq\emptyset$.

This result is very useful. The Kusuoka's model we presented in Example \ref{ExKusuoka} is arbitrage free, because there exists an equivalent change of measure such that $\tau$ is independent from $\F_T$, and hence immersion holds (see \cite{Kusuok} page 79-80 for details). Also, one can show that the $\F_\infty$-measurable random times which are not stopping times do not fulfill this non arbitrage condition.

\begin{lem}
Let $\tau$ be a random time which is $\F_\infty$-measurable. Then, $\mathcal{I}(\P)\neq\emptyset$ if and
 only if $\tau$ is an $\FF$ stopping time (in this case  $\GG=\FF$).
\end{lem}
\proof Suppose that $\exists$ $\P^* \in \mathcal{I}(\P)$. Then $\forall t\geq 0$, $\P^*(\tau>t|\F_t)=\P^*(\tau>t|\F_\infty)$. Now, since $\tau$ is $\F_\infty$ measurable, we have $\P^*(\tau>t|\F_\infty)=\I_{\tau>t}$, and hence $\P^*(\tau>t|\F_t)=\I_{\tau>t}$. This is possible if and only if $\{\tau>t\} \in \F_t$ $\forall t$, that is  if and only if $\tau$ is an $\FF$ stopping time. The converse is obvious. \finproof

\begin{rem}This shows that honest times (which are ends of predictable sets) are not suitable for modeling default events in an arbitrage free financial market of the type introduced in Section \ref{sec::finmodel}. They are encountered in models with insider information, where insiders are shown to obtain free lunches with vanishing risks (\cite{Imkeller02}). Another example of $\F_\infty$-measurable times appears in the models with delayed information.
\end{rem}

Now, we would like to answer the following question: are there more general changes of probability measures that preserve the immersion property? More generally, how is the predictable compensator of $\tau$ modified under an equivalent change of probability measure?  Indeed, it is known that the market implied default intensities (i.e., risk-neutral) are very different from the ones computed using historical data from defaults (i.e., under the physical measure). Hence, for the financial applications it is important to understand how the predictable compensator is modified under general changes of probability measures. Note also the recent paper \cite{ElJeJi} where a particular case is studied: the $\FF$-conditional distribution of $\tau$ admits a density with respect to some non atomic positive measure.

For sake of completeness, we state a general result  due to Jeulin and Yor \cite{yorjeulinens}  which is unfortunately not easy to use in practice:

\begin{prop}[Jeulin-Yor \cite{yorjeulinens}]\label{proplL}
Let $\Q$ be a probability measure which is equivalent to $\P$ on
$\GG$, with $\rho= d\mathbb{Q}/d\mathbb{P}$ on  $\mathcal {G}_\infty$. Define $E$ and $e$ as
in (\ref{densities}) and suppose  that $\FF \overset{\P}{\hookrightarrow} \GG$. Then, $\FF \overset{\Q}{\hookrightarrow}
\GG$ if and only if:
\begin{equation}
\forall t\geq0,\;X\in\mathcal{F}_{\infty},\quad \dfrac{\mathbb{E}_{\mathbf{P}}\left[X\rho|\mathcal{G}_{t}\right]}{E_{t}}=\dfrac{\mathbf{E}_{\mathbf{P}}\left[X\rho|\mathcal{F}_{t}\right]}{e_{t}}.\label{JY_cond}
\end{equation}
In particular,  if
$\rho$ is $\mathcal{F}_{\infty}$-measurable, then $e =E $ and $\FF
\overset{\Q}{\hookrightarrow} \GG$.
\end{prop}

\proof Using Bayes formula, (\ref{JY_cond}) is equivalent to:
$$ \forall t\geq0,\;X\in\mathcal{F}_{\infty},\quad  \mathbb{E}_{\mathbf{Q}}\left[X |\mathcal{G}_{t}\right] = \mathbb{E}_{\mathbf{Q}}\left[X
|\mathcal{F}_{t}\right],$$ which is equivalent to the immersion
property under the measure $\Q$ from Theorem
\ref{lecasdregeneral}. \finproof

\begin{rem}
This theorem holds for more general filtrations (i.e. $\GG$ does not necessarily have to be obtained by progressively enlarging $\FF$ with a random time).  Moreover, although it is not mentioned in \cite{yorjeulinens}, the necessary and sufficient condition is valid even if $\FF$ is not immersed into $\GG$ under $\P$. However, it will not directly help us find a larger class than the change of probability measures for which the density $\rho$ is $\F_\infty$-measurable.
\end{rem}

\section{Some martingale representation properties}

In the remainder of this paper, we suppose that $\tau$ is such that condition \textbf{(A)} holds and that the immersion property holds under $\P$. Recall from Section \ref{sec::tpsaleatoires}, that these assumptions imply that the Az\'ema supermartingale $(Z^\tau_t)$ is a decreasing and continuous process.

Under these assumptions, we prove in this section several general martingale representation theorems for martingales of the larger filtration $\GG$. These results will allow us to construct in section \ref{sec::further} yet larger classes of equivalent probability measures that preserve the immersion property.

We begin with a few useful lemmas.

\begin{lem} \label{projLemma} Assume that \textbf{(A)} and  $\FF \overset{\P}{\hookrightarrow} \GG$ hold. Let  $H$ be a   $\GG$-predictable process and let $N_t=\I_{\tau\leq t}-\Gamma_{t\wedge\tau}$ ($N$ is a $\GG$ martingale). If $\E[|H_\tau|]<\infty$, then:

\begin{equation}
\E\left[\int_0^t H_s dN_s|\F_t\right]=0
\end{equation}
\end{lem}

\proof First we note that if $\E[|H_\tau|]<\infty$, then the
integral $\int_0^\infty |H_s| dN_s$ is well defined. It is enough to check that both integrals
$\int_0^\infty |H_s| d\I_{\tau\leq s}$ and $\int_0^\tau |H_s|
\frac{dA_s^\tau}{Z^\tau_s}$ are finite. The first integral is equal to
$|H_\tau|$ and is hence finite. For the second integral, using the
fact that $A^\tau$ is continuous and hence predictable and using
properties of predictable projections, we have: 
$$\E\left [\int_0^\tau |H_s| \frac{dA_s^\tau}{Z^\tau_{s}}\right]=\E\left[\int_0^\infty \I_{\tau>s}|H_s| \frac{dAs^\tau}{Z^\tau_{s}}\right]=\E\left[\int_0^\infty \;^p(\I_{\tau>s}|H_s|) \frac{dA_s^\tau}{Z^\tau_{s}}\right],$$where $\;^p(\cdot)$
denotes the $(\FF,\P)$-predictable projection. Now, we use the fact that on the interval $s\leq\tau$, $H$ is equal to an $\FF$ predictable process and that $\;^p(\I_{\tau>s})=Z^\tau_{s}$ (because $\tau$ avoids $\FF$ stopping times) to conclude that $\E[\int_0^\tau |H_s| \frac{dA_s^\tau}{Z^\tau_{s}}]=\E[|H_\tau|]$ and consequently the integral $\int_0^\tau |H_s| \frac{dA^\tau_s}{Z^\tau_{s}}$ is also finite.

Since $N$ is a local martingale of finite variation, it is purely
discontinuous. Now, let $(M_t)$ be any square integrable
$\FF$-martingale.  Since $\FF \overset{\P}{\hookrightarrow} \GG$, $(M_t)$ is also a $\GG$-martingale. We also have $[M,N]_t=0$ because $N$ is purely discontinuous, and has a single jump at $\tau$ which
avoids $\FF$ stopping times. Consequently, $N$ is strongly orthogonal to all  $\FF$-martingales, and hence $\E(M_tN_t)=0$ for  all $t$ and all square integrable $\FF$-martingales. This proves the lemma. \finproof

\begin{lem}
[\cite{bremaudyor}] \label{projLemma2} Assume that $\FF
\overset{\P}{\hookrightarrow} \GG$. Let  $H$ be a  bounded
$\GG$-predictable process and let $m$ be an $\FF$ local
martingale. Then:

\begin{equation}
\E^\P\left[\int_0^t H_s dm_s|\F_t\right]=\int_0^t  \;^{(p,\P)}H_s dm_s,
\end{equation}
where $\;^{(p,\P)}H$ is the $(\FF,\P)$-predictable projection of the process $H$.

\end{lem}

We now deduce easily from lemma \ref{projLemma} the following projection formula: 

\begin{lem}\label{lemmedeprojection}Let $\tau$ be any random time.
\begin{enumerate}[(i)]
\item Assume that \textbf{(A)} holds. Then:
$$\E[z_\tau\;\I_{\tau>t}|\F_t]=\E\left[\int_t^\infty z_sdA^\tau_s|\F_t\right]$$
\item Assume further that $\FF \overset{\P}{\hookrightarrow} \GG$. Let $z$ be an $\FF$-predictable process, such that $\E[|z_{\tau}|]<\infty$. Then:
\begin{equation}
\E[z_\tau|\F_t]=\E\left[\int_0^\infty z_sdA^\tau_s|\F_t\right];
\end{equation}if moreover the hazard process $\Gamma$ is defined for all $t\geq0$, that is if $Z_t^\tau>0$ for all $t\geq0$, then
$$\E[z_\tau|\F_t]= \E\left[\int_0^\infty z_s e^{-\Gamma_s}d\Gamma_s|\F_t\right].$$
\end{enumerate}
\end{lem}
\proof (i)This is a consequence of the projection formulae T25, p.
104 in \cite{cdellacherie}; see also \cite{Ashkanstopandnonstop}.

(ii)It is enough to check the result for $z_s=H_r\I_{(r,u]}(s)$, with $r< u$ and  $H_r$  an integrable $\F_r$ measurable random variable. But in this case the result is an immediate consequence of Theorem \ref{lecasdregeneral}.

\finproof

We now state and prove a first representation theorem result for some $\GG$ martingales under the assumption that $(Z_t^\tau)$ is continuous and decreasing, that is  $\tau$ is a pseudo-stopping time that avoids stopping times (the pseudo-stopping time assumption is an extension of the $(H)$ hypothesis framework, see \cite{AshkanYor} and \cite{coculescunikeghbali}). This result was in \cite{bsjeanblanc} (without the pseudo-stopping times there). We give here a simpler proof which easily extends to any random time. But before we state a lemma which we shall use in the proof.

\begin{lem}\label{cor::NTLT}[\cite{bsjeanblanc},  \cite{jenbrutk1}]
Let $\tau$ be an arbitrary random time. Define $$L_t=\I_{t<\tau}e^{\Gamma_t}.$$Then $(L_t)_{t\geq0}$ is a $\GG$ martingale, which is well defined for all $t\geq0$.

Let $\tau$ be a pseudo-stopping time and assume that \textbf{(A)} holds (or equivalently assume that $(Z_t^\tau)$ is continuous and decreasing). Then
$$L_t=1-\int_0^t \frac{\dd N_s}{Z_s^\tau},$$where $(N_t)$ is the $\GG$ martingale $N_t=\I_{\tau\leq t}-\Gamma_{t\wedge\tau}$.
\end{lem}

\begin{thm}\label{premiererepresentation}
Let $\tau$ be a pseudo-stopping time and assume  that \textbf{(A)} holds. Let $z$ be an $\FF$-predictable process such that $\E[|z_\tau|]<\infty$. Then
\[
\E[z_\tau|\G_t]=m_0+\int_0^{t\wedge \tau}\frac{d
m_s}{Z^\tau_s}+\int_0^t (z_s-h_{s})dN_s,
\]
where $ m_t=\E[\int_0^\infty z_s dA^\tau_s|\F_t]$ and $h_t=(Z^\tau_t)^{-1}\left(m_t-\int_0^tz_s dA^\tau_s\right)$.
\end{thm}
\proof It is well known  that (see \cite{delmaismey},
\cite{jenbrutk1} lemma 3.2 or \cite{Ashkansurvey} p.58):
$$\E[z_\tau |\G_t]=L_t\E[z_\tau \I_{\tau>t}|\F_t]+z_\tau \I_{\tau\leq t}.$$ Furthermore, from Lemma \ref{lemmedeprojection}, with the notation of the Theorem,  we have:
$$L_t\E[z_\tau \I_{\tau>t}|\F_t]=m_t-\int_0^tz_s dA^\tau_s.$$ Consequently,
$$\E[z_\tau |\G_t]=L_t m_t-L_t\int_0^tz_s dA^\tau_s +z_\tau \I_{\tau\leq t}.$$Now, noting that $L_t$ is a purely discontinuous martingale with a single jump at $\tau$, we obtain that $L_t$ is orthogonal to any $\FF$ martingale. An integration by parts combined with lemma \ref{cor::NTLT} yields the desired result.
\finproof

\begin{rem}
The proof of Theorem \ref{premiererepresentation} can be adapted so that the result would hold for an arbitrary random time that avoids stopping times. The only thing to modify is lemma \ref{cor::NTLT}: for an arbitrary random time $\tau$ that avoids stopping times, $Z_t^\tau$ is continuous and not of finite variation anymore, so that an extra term must be added when expressing $L_t$ as a sum of stochastic integrals.
\end{rem}

Now we state a corollary which will play an important role in our search for a larger class of equivalent probability measures which preserve the immersion property.

\begin{cor}\label{projsurgt}
Let $\tau$ be a random time such that \textbf{(A)} and $\FF \overset{\P}{\hookrightarrow} \GG$ hold. Let $z$ be an $\FF$ predictable process such that $\E[|z_\tau|]<\infty$. Assume further that  there exists a constant $c$ such that $\E[z_\tau|\F_\infty]=c$. Then, there exists an $\FF$-predictable process $(k_{t})$, such that:
\[
\E[z_\tau|\G_t]=c+\int_0^t k_s dN_s.
\]
\end{cor}

\proof  Using the fact that $Z^{\tau}$ is continuous and Theorem
\ref{premiererepresentation}, we have (we also use the fact that since $\tau$ avoids stopping times, we can replace $h_{s-}$ with $h_s$):

\[
\E[z_\tau|\G_t]=m_0+\int_0^{t\wedge \tau}\frac{d m_s}{Z^\tau_s}+\int_0^t (z_s-h_s)dN_s,
\]
where $ m_t=\E[\int_0^\infty z_s dA^\tau_s|\F_t]$ and $h_t=(Z^\tau_t)^{-1}\left(m_t-\int_0^tz_s dA^\tau_s\right)$.
Now, from lemma \ref{lemmedeprojection}, (ii), we also have under the assumptions of the corollary that
$$m_t=\E[\int_0^\infty z_s dA^\tau_s|\F_t]=\E[z_\tau|\F_t].$$ Since it is assumed that $\E[z_\tau|\F_t]=\E[\E[z_\tau|\F_\infty|\F_t]]=c$, the result of the corollary follows at once, with $k_t=z_t-h_t$.
\finproof 

We now combine corollary \ref{projsurgt} with
Proposition \ref{invariancecomp} to obtain a representation
theorem for a larger class of $\GG$ martingales.

\begin{prop}\label{genrepresentation}
Let $\tau$ be a random time such that \textbf{(A)} and $\FF
\overset{\P}{\hookrightarrow} \GG$ hold. Let $G=Fz_\tau$, where
$F$ is an integrable, $\F_\infty$-measurable random variable such
that $F\neq0,\; \text{ a.s. }$ and $z$ is an $\FF$-predictable
process, such that $z_\tau F$ is integrable. Then:
\begin{equation}
\E[G|\G_t]=\E[G]+\int_0^t \left(\E[G]+ Y_{s} - L_s
\frac{m^G_s}{m^F_s}+\int_0^s k_u d N_u\right)dm^F_s +\int_0^{t}
L_s dm^{G}_s+\int_0^t m^F_s k_s d N_s,
\end{equation}
where:
\[
m^F_t:=\E[F|\F_t]; \quad m^G_t:=\E[G|\F_t]; \quad Y_{t }=\int_0^{t} L_s d\left(\frac{m^{G}_t}{m^F_t}\right),
\]and  where $(k_t)$ is an $\FF$-predictable process (which can be given explicitly).
\end{prop}

\proof Without loss of generality, we can assume  that $F$ is
strictly positive and that  $\E[F]=1$ (the general case would
follow by writing $F=F^+-F^-$). Then, we define  $d\tilde{\Q}\vert_{\G\infty}=F \cdot d\P\vert_{\G\infty}$.

 Hence, from Proposition \ref{invariancecomp}, the $(H)$
hypothesis holds under $\tilde{\Q}$ and $\tilde{\Q}(\tau >
t|\F_t)=\P(\tau >t|\F_t)$. We then obtain:

\[
\E[G|\G_t]=\E[z_\tau F|\G_t]=\E[F|\G_t]\E^{\tilde{\Q}}[z_\tau| \G_t]=m^F_t\E^{\tilde{\Q}}[z_\tau|\G_t],
\]
Using the decomposition from Theorem \ref{premiererepresentation}, we get:
\[
\E^{\tilde\Q}[z _\tau|\G_t]=\E^{\tilde \Q}[z_\tau]+Y_{t}+\int_0^t k_s dN_s,
\]
where $Y_{t }=\int_0^{t } L_s d  \widetilde{m}_s$. Here,  $\tilde
m$ is the $\tilde \Q$-martingale defined by
\[
\tilde m_t:=\E^{\tilde\Q}[z(\tau)|\F_t]=\E^{\P}[z(\tau)F|\F_t](m^F_t)^{-1}=\frac{m^{G}_t}{m^F_t}.
\]
 and
$$k_t=z_t-(Z^\tau_t)^{-1}\left(\widetilde{m}_t-\int_0^tz_s dA^\tau_s\right).$$
Consequently:

$$\E[G|\G_t] =m^F_t \left( \E^\P[G]+\int_0^{t} L_s d\left(\frac{m^{G}_s}{m^F_s}\right)+\int_0^t k_s dN_s\right).$$
Now, an integration by parts formula and  some tedious computation lead to:

$$\E[G|\G_t]=\E[G]m^F_t+\int_0^t \left(
Y_{s} - L_s\frac{m^{G}_s}{m^F_s}+\int_0^s k_ud N_u\right)dm^F_s
+\int_0^{t} L_sdm^{G}_s+\int_0^t m^F_s k_sd N_s,$$which completes
the proof of our theorem. \finproof

As a corollary, we obtain the following generalization of a representation result by Kusuoka \cite{Kusuok}, which was obtained in the Brownian filtration.

\begin{cor}\label{corgenrepresentation}
Let $\tau$ be a random time such that \textbf{(A)} and $\FF \overset{\P}{\hookrightarrow} \GG$ hold. Then any $\GG$-locally square integrable martingale $(M_t)$ can be written as:
\begin{equation}\label{orth_decomp}
M_t=M_0+V_t+\int_{0}^{t}h_s \dd N_s,
\end{equation}
where $(V_t)$ is in the closed subspace of $\GG$-locally square integrable martingales generated by the stochastic integrals of the form $\int_0^t R_s \dd m_s$, where $(m_t)$ is an $\FF$-locally square integrable martingale, $(R_t)$ is a $\GG$-predictable process such that $\int_0^t R_s^2 d\langle m,m\rangle_s$ is locally integrable, and where $(h_t)$ is an $\FF$-predictable process which is such that $h_\tau^2$ is  integrable.
\end{cor}
\proof  The result follows from Proposition
\ref{genrepresentation} and the fact that any
$\G_\infty$-measurable random variable can be written as a limit
of finite linear combinations of functions of the form $Ff(\tau)$
where $F$ is an $\F_\infty$ random variable and $f$ a Borel
function such that $Ff(\tau)$ is integrable. \finproof
\begin{rem}
Since any element $V$ in the closed subspace of $\GG$-locally square integrable martingales generated by the stochastic integrals of the form $\int_0^t R_s \dd m_s$ is strongly orthogonal to the purely discontinuous martingales of the form $\int_{0}^{t}h_s \dd N_s$, it follows that the  decomposition (\ref{orth_decomp}) is unique.
\end{rem}
\begin{cor}\label{corfiltrationbrownienne}\cite{Kusuok}
Assume that $\FF$ is the natural filtration of a one dimensional Brownian motion $(W_t)$. Let $\tau$ be a random time such that \textbf{(A)} and $\FF \overset{\P}{\hookrightarrow} \GG$ hold. Then any $\GG$-locally square integrable martingale $M$ can be written as:
$$M_t=M_0+\int_{0}^t R_s \dd W_s+\int_{0}^{t}h_s \dd N_s,$$
where $(R_t)$ is a $\GG$-predictable process such that $\int_0^t R_s^2 \dd s$ is locally integrable, and where $(h_t)$ is an $\FF$-predictable process which is such that $h_\tau^2$ is  integrable.
\end{cor}
\begin{rem}
A  result similar to the representation of corollary \ref{corfiltrationbrownienne} would hold if the filtration $\FF$ has the predictable representation property with respect to a family of locally square integrable martingales.
\end{rem}

Combining Lemma \ref{projLemma} and Corollary \ref{corgenrepresentation} one gets:
\begin{cor}\label{correprorthf} Let $\tau$ be a random time such that \textbf{(A)} and $\FF \overset{\P}{\hookrightarrow} \GG$ hold. Assume that $(M_t)$ is a locally $L^2$  $\GG$-martingale. $(M_t)$ is strongly orthogonal to all locally $L^2$ $\FF$-martingales if and only if there exists  an $\FF$-predictable process $(h_t)$, such that $h_\tau^2$ is  integrable and  such that
\[
M_t=M_0+\int_{0}^{t}h_s \dd N_s.
\]
\end{cor}

\section{Equivalent changes of probability measures: further results}\label{sec::further}

In this section, we prove two important results. We first characterize the Radon-Nikod\'ym derivative $\dd\mathbf{Q}/\dd\mathbf{P}$ of the measures $\Q \in \mathcal I(\P)$. Then, we generalize Proposition \ref{invariancecomp}: we compute the Az\'ema supermartingale of a random time for which the immersion property holds for a very large class of equivalent change of probability measures.

We begin with  a lemma which is of interest for its own sake:

\begin{lem}\label{lemIntermed}
Let $$\dd\mathbf{Q}/\dd\mathbf{P}=H \quad \text{on }\G_\infty,$$ where $H$ is a positive and integrable $\F_\tau$-measurable random variable such that $\E^\P[H|\F_\infty]=1$. Then $\FF\overset{\Q}{\hookrightarrow} \GG$ holds.
\end{lem}
\begin{proof}
From Corollary \ref{projsurgt} we have that $E_t:=\E^\P[H|\G_t]= 1+\int_0^t h_s dN_s$ and  $\E^\P[H|\F_t]=1$. In addition, since $\tau$ avoids the $\FF$-stopping times and since $E$ is a purely discontinuous martingale, $[M,E]=0$, for any $(\FF,\P)$-martingale $(M_t)$.
Hence, by  Girsanov's theorem (H) holds under $\Q$.
\end{proof}

\begin{thm}
Let $\tau$ be a random time such that \textbf{(A)} and $\FF\overset{\P}{\hookrightarrow} \GG$ hold. Let $\Q$ be a probability measure which is equivalent to $\P$.
\begin{enumerate}[(i)]
\item Assume that 
$$\dd\mathbf{Q}/\dd\mathbf{P}=FH \quad \text{on }\G_\infty,$$where $F$ is  a positive  $\mathcal{F}_{\infty}$-measurable and integrable random variable with $\E[F]=1$, and $H$ is a positive   $\F_\tau$ measurable and integrable random variable such that $\E^\P[H|\F_\infty]=1$ (and such that $FH$ is integrable). Then $\FF\overset{\Q}{\hookrightarrow} \GG$ holds.
\item Conversely, assume that $\FF\overset{\Q}{\hookrightarrow} \GG$ holds. With the notation (\ref{densities}), assume further that 
\begin{equation}\label{condL2}
\E\left[\dfrac{e_\infty^2}{E_\infty}\right]<\infty.
\end{equation}
Then, there exist $F$ and $H$ as in $(i)$ above, such that 
$$\dd\mathbf{Q}/\dd\mathbf{P}=FH \quad \text{on }\G_\infty.$$
\end{enumerate}
\end{thm}

\begin{proof}
$(i)$ Assume that $H$ is $\F_\tau$-measurable.  Introduce: $\dd \tilde{\Q} =F \cdot \dd\P$, hence $\dd \Q=H\cdot  d\tilde{\Q}$ and notice that $\E^{\tilde Q}[H|\F_\infty]=1$.  From Proposition \ref{invariancecomp}, we know that $(H)$ holds under $\tilde{\Q}$, then using Lemma \ref{lemIntermed}, it follows that the immersion property also holds under $\Q$.

$(ii)$ Recall first the following general fact from Theorem \ref{JY_decomp}: if $\FF\overset{\P}{\hookrightarrow} \GG$ holds, then $\eta_t:=e_t/E_t$ is a $(\GG,\Q)$ uniformly integrable martingale. We then note that:
$$\E^\P[(\eta_\infty)^{-1}|\F_\infty]=\E^\Q[(\eta_\infty)^{-1} (E_\infty)^{-1}|\F_\infty]e_\infty=1.$$
Since $E_t=e_t( \eta_t)^{-1}$, it follows that  $\dd\mathbf{Q}/\dd\mathbf{P}=E_\infty= FH$, where $F=e_\infty$ is $\FF$-measurable with $\E[F]=1$ and $H=(\eta_\infty)^{-1}$ satisfies $\E^\P[H|\F_\infty]=1$.

Now, let us assume further that $\FF\overset{\Q}{\hookrightarrow} \GG$ also holds. Assumption (\ref{condL2}) is easily seen to mean that $(\eta_t)$ is an $L^2(\GG,\Q)$ bounded martingale. Using twice Girsanov's theorem, one can show that if $(m_t)$ is any $L^2(\FF,\Q)$ bounded martingale, then $(m_t\eta_t)$ is a  $(\GG,\Q)$ uniformly integrable martingale. Indeed, if $(m_t)$ is an $(\FF,\Q)$  martingale, then, from Girsanov's theorem, $(m_t e_t)$ is an $(\FF,\P)$ martingale. Now, because $\FF\overset{\P}{\hookrightarrow} \GG$ holds, we also have that $(m_t e_t)$ is an $(\GG,\P)$. Now another application of Girsanov's theorem yields that $m_t\frac{e_t}{E_t}$, which is (by definition) $(m_t\eta_t)$, is a $(\GG,\Q)$ martingale. In other words, the $(\GG,\Q)$ martingale $\eta$ is strongly orthogonal to all $(\FF,\Q)$-martingales viewed as $(\GG,\Q)$ martingales (recall that by assumption $\FF\overset{\Q}{\hookrightarrow} \GG$). Then, by Corollary  \ref{correprorthf}, $\eta_\infty$ is $\F_\tau$-measurable, and so is $H=(\eta_\infty)^{-1}$. 
\end{proof}

The Az\'ema supermartingale plays an important role in credit risk modeling. Now, we would like to display the form of the  $\Q$-Az\'ema supermartingale, denoted $Z^\Q$, under a large class of  equivalent change of probability measures. Before doing so, we would like to state a very useful, though somehow forgotten, result by It\^o and Watanabe \cite{itowatanbe} on multiplicative decompositions of supermartingales. In particular, the multiplicative decomposition reveals to be useful in the study of the intensity of the default time as we shall see.

\begin{thm}[It\^{o}-Watanabe \cite{itowatanbe}] Let  $\left(Z_{t}\right)$ be a nonnegative c\`{a}dl\`{a}g
supermartingale, and define
$$T_{0}=\inf\left\{t:\;Z_{t}=0\right\}.$$
Suppose $\mathbb{P}\left(T_{0}>0\right)=1$. Then $Z$ can be factorized as:$$Z_{t}=Z_{t}^{(0)}Z_{t}^{(1)},$$with a positive local martingale $Z_{t}^{(0)}$ and a decreasing process $Z_{t}^{(1)}$ ($Z_{0}^{(1)}=1$). If there are two such factorizations, then they are identical in $[0,T_{0}[$.
\end{thm}

It follows that, if $\forall t$ $Z^\tau_t>0$ $a.s.$, and is continuous, then there exist a unique local martingale $(m^\tau_t)$ and a unique predictable increasing process $(\Lambda_t)$ such that: $$Z^\tau_t=\mathcal{E}(\int_0^\cdot \frac{dm_s^\tau}{Z_s})_t e^{-\Lambda_t},$$ where the process $\Lambda$ is given by: $$\Lambda_t=\int_{0}^t\frac{1}{Z_{s}^{\tau}}da_{s}^{\tau}.$$ From Theorem \ref{calccomp} the process: $$N_t:=\I_{\{\tau\leq t\}}-\Lambda_{t\wedge\tau}$$ is a $\GG$ martingale.

\begin{thm}\label{PAzema} Assume  that $\FF \overset{\P}{\hookrightarrow} \GG$ and that $\tau$ is a random time that avoids stopping times. Assume further that $Z^\tau_t>0$ for all $t\geq0$.
Let $(m_t)$ be an $(\FF,\P)$-martingale and let $F$ be a $\GG$ predictable processes such that $\mathcal{E}\left(\int_0^\cdot F_s dm_s \right)_t$ is a uniformly integrable $\GG$-martingale. Let $H$ be an $\FF$ predictable process such that $\mathcal{E}\left(\int_0^\cdot H_s dN_s \right)_t$ is uniformly integrable $\GG$-martingale. Let
\[
E_t= \mathcal{E}\left(\int_0^\cdot F_s dm_s \right)_t\mathcal{E}\left(\int_0^\cdot H_s dN_s \right)_t \label{density}.
\]  Assume further that $(E_t)$ is a uniformly integrable $\GG$-martingale (this is the case for example if $\mathcal{E}\left(\int_0^\cdot F_s dm_s \right)_t$ and $\mathcal{E}\left(\int_0^\cdot H_s dN_s \right)_t$ are bounded in $L^2$, since $\int_0^t F_s dm_s $ and $\int_0^t H_s dN_s$ are orthogonal). Define $$\dd\Q=E_t\cdot \dd\P \text{ on }\G_t.$$

Then, the $\Q$-Az\'ema supermartingale associated with $\tau$ has the following multiplicative decomposition:
\[
Z^\Q_t=\Q(\tau>t|\F_t)=\mathcal{E}\left(\int_0^\cdot \left( \tilde{F}_s -\;^{(\Q, p)}F_s\right)d\widetilde{m}_s\right)_t e^{-\int_0^t(1+H_s)d\Gamma_s}
\]
where:

\begin{itemize}
\item $^{(\Q,p)}F$ is the $\FF$-predictable projection of the process $F$ under the probability $\Q$;
\item $\tilde{F}$ is an $\FF$-predictable process such that $\I_{\{\tau > t \}}F_t=\I_{\{\tau > t \}}\tilde F_t $ and
\item $\widetilde{m}_t=m_t-\int_{0}^{t}\frac{\dd[m,e]_s}{e_{s-}}$ is a $(\Q,\FF)$-martingale.
\end{itemize}

It follows that the process:
\[
N^\Q_t:= \I_{\{\tau\leq t\}}-\int_0^{t\wedge\tau}(1+ H_s)d\Lambda_s
\]
is a $(\GG,\Q)$-martingale.
In particular, if the process $F$ is $\FF$-predictable, then:
\[
Z^\Q_t=\Q(\tau>t|\F_t)=e^{-\int_0^t(1+H_s)d\Lambda_s},
\]
and the immersion property holds under $\Q$.
\end{thm}
\begin{rem}
The process $H$ above is taken to be $\FF$ measurable to simplify the notations. Indeed, since the martingale $N$ is constant after $\tau$, and since a $\GG$ predictable process before $\tau$ is equal to an $\FF$ predictable process, we could as well take $H$ to be $\GG$ predictable.
\end{rem}
\begin{proof}
First, we need to compute $e_t:=\E^\P[E_t|\F_t]$. When applying the Lemma \ref{projLemma} to:

\[
E_t=1+ \int_0^t E_{s-} F_s dm_s +\int_0^t E_{s-} H_s dN_s
\]
we obtain that:

\[
e_t=1+ \int_0^t  \;^{(p,\P)}(E_{s-} F_s) dm_s
\]
with $\;^{(p,\P)}(E_{t-} F_t)=\E^\P[E_{t-}F_t|\F_{t-}]=\E^\Q[F_t|\F_{t-}]e_{t-}$. Hence: $$e_t= \mathcal{E} (\int_0^t\;^{(p,\Q)}F_s dm_s).$$

Replacing this  in the formula: $Z^\Q_t=\E^\P[\I_{\{\tau > t\}}E_t|\F_t]/e_t$ leads us to: 

$$Z^\Q_t =e^{-\int_0^t(1+H_s)d\Lambda_s}  \frac{\mathcal{E}(\int_0^t \widetilde F_s dm_s) }{\mathcal{E} (\int_0^t \;^{(p,\Q)}F_s dm_s)}=\exp\left\{ \int_0^t ( \widetilde F_s -  \;^{(p,\Q)}F_s )dm_s-\frac{1}{2}\int_0^t   ( \widetilde F_s ^2-  \;(^{(p,\Q)}F_s) ^2)d[m,m]_s\right\}.$$
Using Girsanov's theorem, $\widetilde{m}_t=m_t-\int_{0}^{t}\frac{\dd[m,e]_s}{e_s-}=m_t-\int_0^t \;^{(p,\Q)}F_s d[m,m]_s$ is an $(\FF,\Q)$-martingale. The result follows when replacing  $m_t$ in the above expression of $Z^\Q_t$ by $\widetilde{m}_t+\int_0^t \;^{(p,\Q)}F_s d[m,m]_s$.
\end{proof}
\begin{cor}
Suppose that $\FF \overset{\P}{\hookrightarrow} \GG$ and that  \textbf{(A)} hold. Assume further that $Z_t>0$ for all $t\geq0$. Define $\Q$ on $\G_t$ by:
\[
\dd\Q/\dd\P=E_t:= \mathcal{E}\left(\int_0^\cdot f_s dm_s \right)_t,
\]
with $f$ a  $\GG$-predictable process such that $E$ is a uniformly integrable martingale. Then, under $\Q$, the process $N_t=\I_{\{\tau \leq t \}}- \Gamma_{t \wedge \tau}$ remains a $\GG$-martingale.
\end{cor}
\begin{proof}
It suffices to take $H=0$ in theorem \ref{PAzema}.
\end{proof}

\renewcommand{\refname}{References}

\end{document}